# G.C.H. implies existence of many rigid almost free abelian groups[*]

Rüdiger Göbel
Fachbereich 6
Mathematik
Universität Essen
45117 Essen, Germany

e-mail:
R.Goebel@uni-essen.de

Saharon Shelah
Department of Mathematics
Hebrew University
Jerusalem, Israel
and Rutgers University
New Brunswick, N. Y., U.S.A.

e-mail:
Shelah@math.huji.ae.i1

October 28, 2018

## §1 Introduction

More than thirty years ago A. L. S. Corner [2] proved a very remarkable theorem concerning the existence of "pathological" abelian groups; it is formulated in terms of endomorphism rings.

*If $A$ is a countable ring with reduced and torsion-free additive group $A_{\mathbb{Z}}$, then there exists a reduced and countable abelian group $G$ with endomorphism ring $End_{\mathbb{Z}} G = A$.*

This result has been generalized in various directions, e.g. replacing the ground ring $\mathbb{Z}$ by a more general ring and dropping the cardinal restriction (in Corner's result this is actually $|A| < 2^{\aleph_0}$), see [2], [20], [8], [9], [4], [23].

Other extensions include abelian groups which are not necessarily torsion-free, see Corner, Göbel [4] for references.

Here we are interested in a different strengthening of Corner's important result: As a byproduct of older constructions of abelian groups $G$ with prescribed endomorphism ring $A$ we have obtained (in a strong sense) almost-free groups provided $A_{\mathbb{Z}}$ is free and provided we are working in Gödel's universe $V = L$.

[*]This work is supported by GIF project No. G-0294-081.06/93 of the German-Israeli Foundation for Scientific Research & Development

Let us recall the appropriate definition. An $R$-module $G$ (of cardinality $\lambda$) with $|R| < |G|$ is strongly $\lambda$-free if $G$ has a $\lambda$-filtration (a continuous chain of submodules $\{G_\nu : \nu < \lambda\}$ with $|G_\nu| < \lambda$ and $\bigcup_{\nu < \lambda} G_\nu = G$) with two additional properties:

(a) Any $G_\nu$ is a free $R$-module.

(b) If $\nu$ is a successor ordinal and $\nu < \rho < \lambda$, then $G_\rho/G_\nu$ is free as well.

We say that $G = \bigcup_{\nu < \lambda} G_\nu$ is a $\lambda$-filtration for strong $\lambda$-freeness. First results on almost-free groups of cardinality $2^{\aleph_0}$ are in [1].

Because strong $\lambda$-freeness is a very natural definition, it can be reformulated in many ways, see Eklof, Mekler [15] for equivalent characterizations and many results on this class of modules. In order to illustrate such results we will restrict ourselves to abelian groups, i.e. $R = \mathbb{Z}$. In $V = L$ we have the following strengthening of Corner's result indicated above.

*If $A$ is a ring, $A_\mathbb{Z}$ is a free abelian group and $\lambda > |A|$ is any regular, not weakly compact cardinal, then we can find a strongly $\lambda$-free abelian group $G$ with $|G| = \lambda$ and $End_\mathbb{Z} G = A$.*

This older result [8] grew out of ideas of Eklof, Mekler [14], and Shelah [38]. It is closely related to the work of Jensen [30]. It is only natural to ask if such almost free groups exist without the vehicle of additional set theory. Extending a result of Shelah [36], Eda [11] showed that there is at least some hope.

There exists an $\aleph_1$-free abelian group $G$ of cardinality $\aleph_1$ with trivial dual $G^* = 0$, see also Corner, Göbel [5]. Recall that $Hom(G, R) = G^*$ is the dual of an $R$-module, and $G$ is $\kappa$-free for some cardinal $\kappa$ if any submodule $U$ with $|U| < \kappa$ is contained in a free submodule. In [22] we showed the following stronger result.

*If $A$ is a countable ring and $A_\mathbb{Z}$ is free, then we can find an $\aleph_1$-free abelian group $G$ of cardinality $\aleph_1$ with endomorphism ring $End_\mathbb{Z} G = A$.*

This result provides a satisfying answer in ZFC, but unfortunately (from the algebraist's point of view) it remains restricted to the world $\aleph_1$.

In fact, assuming ZFC + Martin's axiom, it was shown in [22] that $\aleph_2$-free groups of size $\aleph_2$ are already always separable (= pure subgroups of products $\mathbb{Z}^\mathbb{I}$), hence they must split into various summands. Their endomorphism rings will never be $\mathbb{Z}$ for instance. This brings us to the starting point of this paper.

We want to derive a realization theorem of rings $A$ with $A_\mathbb{Z}$ free as $A = End_\mathbb{Z} G$ for strongly $|G|$-free abelian groups $G$ with $|G| \geq \aleph_1$, using as weak

additional set theoretic conditions as possible. From the result mentioned above we know that Martin's axiom must be excluded. Our results discussed below are surprising in two ways. First of all, we will only need a special case of G.C.H., which ensures the existence of the desired modules. In another paper [43], Shelah and Spasojevic will show by proper forcing arguments that even G.C.H. will not be strong enough for $End\, G = \mathbb{Z}$, $|G| = \lambda$ and $G$ $\lambda$-free if $\lambda$ is too large, e.g. if $\lambda$ is any strongly inaccessible cardinal. Hence we are limited to "small large cardinals". Secondly it turns out that the restriction to G.C.H. forces us to develop novel ideas to get the desired modules. A new Step-Lemma is needed! Here we will also apply combinatorial tools fortunately developed by Shelah [35], [41] a decade ago. The main results of this paper can be summarized as follows.

In order to construct pathological groups we begin with the existence of groups with trivial duals for cardinals $\aleph_n$ $(n \in \omega)$. Then we derive results about strongly $\aleph_n$-free abelian groups of cardinality $\aleph_n$ $(n \in \omega)$ with prescribed free, countable endomorphism ring. Finally we use combinatorial results [35], [41] to give similar answers for cardinals $> \aleph_\omega$. As in Magidor and Shelah [31], a paper concerned with the existence of $\kappa$-free, non-free abelian groups of cardinality $\kappa$, the induction argument breaks down at $\aleph_\omega$. Recall that $\aleph_\omega$ is the first singular cardinal and such groups of cardinality $\aleph_\omega$ do not exist by the well-known Singular Compactness Theorem [34], see also Hill [28] for cardinals cofinal to $\omega$ and, for a survey, Eklof and Mekler [15, p. 107].

We will fix a countable commutative ring $R$ with multiplicatively closed subsets $S$ of regular elements containing 1. An $R$-module $G$ will be called torsion-free (with respect to $S$) if $sg = 0$ $(s \in S,\ g \in G)$ implies $g = 0$. Moreover $G$ is reduced (with respect to $S$) if $\bigcap_{s \in S} sG = 0$. If $A$ is an $R$-algebra, and the $R$-module $A_R$ is $R$-free (torsion-free or reduced) we will say that the algebra $A$ is free (torsion-free or reduced). We will also fix a countable, free $R$-algebra $A$ and derive the following results.

**Theorem 5.1:** *Let $(R, S)$ be as above. If $H$ is a strongly $\mu$-free $R$-module of regular cardinality $\mu$ with trivial dual $H^* = 0$ and $\lambda = 2^\mu = \mu^+$ (the successor cardinal of $\mu$), then there exists a strongly $\lambda$-free $R$-module of cardinality $\lambda$ with trivial dual.*

As a consequence of (5.1) we obtain the existence of such $R$-modules of cardinality $\mu_n$ $(n \in \omega)$ for the $n$-th successor $\mu_n$ of $\mu$ provided $2^{\mu_i} = \mu_{i+1}$ for all $i < n$. The existence of $H$ as in (5.1) for $|H| = \aleph_1$ under $2^{\aleph_0} = \aleph_1$ is well-known, see Dugas [7], Shelah [38] or Dugas, Göbel [8], cf. also Eklof, Mekler [15, p. 391]. Hence the existence of the derived strongly $\aleph_n$-free $R$-modules of cardinality $\aleph_n$ with trivial dual is immediate (see also (5.2)).

In §6 we will work with a stronger **algebraic** hypothesis, which again is satisfied automatically for $\aleph_1 = 2^{\aleph_0}$ by the papers [8],[38] above. We derive the

**Theorem 6.1:** *Let $(R, S)$, $\lambda = 2^\mu = \mu^+$ and $\mu$ regular as above. If $H$ is a strongly $\mu$-free $R$-module of cardinality $\mu$ with $End_R H = A$ for some countable, free $R$-algebra $A$, then there exists a strongly $\lambda$-free $R$-module $G$ of cardinality $\lambda$ with $End_R G = A$.*

Now the existence of strongly $\aleph_n$-free $R$-modules $G$ of cardinality $\aleph_n$ with $End_R G = A$ follows for each countable, free $R$-algebra $A$, provided $\aleph_{i+1} = 2^{\aleph_i}$ for each $i < n$. As one of many consequences we see that under the mild restriction $\aleph_{i+1} = 2^{\aleph_i}$ $(i < n)$ there are counterexamples to Kaplansky's test problems for each cardinal $\aleph_n$ $(n \in \omega)$ even ones which are strongly $\aleph_n$-free (see [3] for the required rings). Similarly, we obtain indecomposable, strongly $\aleph_n$-free groups of cardinality $\aleph_n$. By the result of Shelah, Spasojevic [43] we know that a realization theorem similar to (6.1) and even the existence of modules with trivial dual (5.1) does not follow from G.C.H. for large enough cardinals, e.g. for the first strongly inaccessible cardinal $\xi$. What happens at smaller cardinals, for instance at $\aleph_{\omega+1}$? The Theorems 6.1 and (5.1) can be extended to $\aleph_{\omega+1}$, hence to $\aleph_{\omega+n}$ $(n \in \omega)$ by (5.1) and (6.1), and to certain larger cardinals $< \xi$. We must strengthen the methods from §§3,4 and 5, and this procedure is sketched in §7.

# §2 Representation of strongly $\mu$-free modules and the notion of a type of freeness

We fix some notations throughout this paper.

**(2.1)** *Let $R$ be a commutative ring and $S \subseteq R$ a countable, multiplicatively closed subset of regular elements such that $R_R$ (as $R$-modules) is reduced and torsion-free (with respect to $S$), compare §1. We will also assume $1 \in S$ and enumerate $S = \{s_n : \; n \in \omega\}$ and set $q_n = \prod_{i=1}^{n} s_i$.*

Our results concern almost free $R$-modules (mainly for $|R| = \aleph_0$), hence we will refine the notion of almost free and strongly $\mu$-free $R$-modules from §1:

We say that a strongly $\mu$-free $R$-module $H$ has freeness type $ft H = u$ if $u = \langle \mu_1, \cdots, \mu_n \rangle$ is a descending sequence of regular cardinals $\mu_i$ with $\mu_1 = \mu$, $\mu_n = \aleph_0$ and if $H = \bigcup_{\alpha < \mu_1} H_\alpha$ is a $\mu_1$-filtration for strong $\mu_1$-freeness with $H_{\alpha+1}/H_\alpha$ either free or of freeness type $\langle \mu_2, \cdots, \mu_n \rangle$. The definition is completed by induction if we say that countable groups which are not free are of freeness type $\langle \aleph_0 \rangle$. If $H$ is of freeness type $u = \langle \mu_1, \cdots, \mu_n \rangle$, then we can choose a particular filtration of $H$, which we will call a $u$-filtration of $H$. We will write $H = \bigcup_{\alpha_1 < \mu_1} H_{\alpha_1}$, such that $H_{\alpha_1+1}/H_{\alpha_1}$ is either free or of freeness type $\langle \mu_2, \cdots, \mu_n \rangle$ and set $\boldsymbol{S}_{\mu_1} = \{\alpha \in \mu_1 : H_{\alpha+1}/H_\alpha$ is not free $\}$. By inserting $< \mu_1$ members between $H_\alpha$ and $H_{\alpha+1}$ we may assume that $H_{\alpha+1}/H_\alpha \cong R$ if $\alpha \notin \boldsymbol{S}_{\mu_1}$. Moreover **we assume** $\boldsymbol{S}_{\mu_1} \subseteq E_{\mu_2} := \{\alpha \in \mu_1, \, cf\alpha = \mu_2\}$. If $\alpha \in \boldsymbol{S}_{\mu_1}$, then $H_{\alpha+1}/H_\alpha$ is not free and we proceed similarly. We can find a continuous chain

$$H_\alpha = H_{\alpha 0} \subseteq ... \subseteq H_{\alpha\beta} \subseteq ... \subseteq H_{\alpha\mu_2} = H_{\alpha+1} \quad (\beta < \mu_2)$$

with quotients $H_{\alpha\beta+1}/H_{\alpha\beta}$ either free or of type $\langle \mu_3, \cdots, \mu_n \rangle$. Recursively we unravel all of $\langle \mu_1, \cdots, \mu_n \rangle$ and let $u^i = \mu_1 \times \cdots \times \mu_i$ $(i = 1, \cdots, n)$. We obtain a refined, continuous chain $\{H_\alpha, \, \alpha \in u^i\}$ of free $R$-modules $H_\alpha$ with $H_\beta/H_{\alpha+1}$ free for $\beta \geq \alpha + 1$. Moreover $H_{\alpha+1}/H_\alpha$ is either $\cong R$ or of type $\langle \mu_{i+1}, \cdots, \mu_n \rangle$. The set $u^i$ is ordered lexicographically and each $\alpha \in u^i$ is of the form $\alpha = \alpha_1^{\wedge} \cdots^{\wedge} \alpha_i$ with $\alpha_j < \mu_j$. As usual, $\alpha + 1$ will denote the successor of $\alpha$ in $u^i$ which is $\alpha + 1 = \alpha_1^{\wedge} \cdots^{\wedge} \alpha_i + 1$, and $0$ denotes the initial element of $u^i$.

Moreover, again refining the chain $\{H_\alpha : \alpha \in u^n\}$, we may assume that $H_{\alpha+1}/H_\alpha \subseteq S^{-1}R$ for each $\alpha \in u^n$ and $H_0 = 0$. We will identify $u^n$ with $\mu$ in §§3 − 6 and write for $\boldsymbol{S} = \{\alpha \in \mu : H_{\alpha+1}/H_\alpha$ is not free $\}$, hence

$$\textbf{(2.2)} \quad \begin{cases} \text{either } H_{\alpha+1} = H_\alpha \oplus t_{\alpha 0}R, \text{if } \alpha \notin \boldsymbol{S} \text{ or } H_{\alpha+1} = \langle H_\alpha, t_{\alpha k} : k \in \omega \rangle \\ \text{is a free R-module with } q_{\alpha k} t_{\alpha k+1} = t_{\alpha k} + \sum_{j<k_\alpha} t_{\alpha j k_j} \text{ where} \\ q_{\alpha k} \in S, \, t_{\alpha j k_j} \in H_{\alpha j+1} \text{ for } \alpha_j < \alpha \text{ and } k < \omega. \end{cases}$$

**Remark:** The reader not interested in particular freeness-type or in modules of rank $\aleph_{\omega+1}$ may choose (2.2) immediately from $\{H_\alpha, \, \alpha \in u^1\}$. Recall that $H$ is uniquely determined (up to isomorphism) by these generators and relations (as follows by a transfinite induction on $\alpha < \mu$). Inductively we will also use particular free resolutions.

**Definition 2.3:** *If $u = \langle \mu_1, \cdots, \mu_n \rangle$ is a (finite) descending sequence of regular cardinals, then we say that the $R$-module $M$ has a $u$-resolution for $\kappa \in u$, if $\mu_1 < \mu$ and the following holds:*

(1) *$M$ is strongly $\mu$-free of cardinality $\mu$.*

(2) *We have a free resolution $M \cong F/K$.*

(3) *$K = \bigoplus_{k \in \kappa} K_k$ and $F, K_k$ are free $R$-modules of rank $\mu$.*

(4) *$K_k$ has a basis $\{x_{\alpha k} : \alpha < \mu\}$ and if $K_{\beta k} = \langle x_{\alpha j} : j < k$ or $\alpha < \beta \rangle$ for $(k < \kappa)$ $(\beta < \mu)$, then $F/K_{\beta k}$ is free.*

Sometimes $R$ in (2.3) will be replaced by an $R$-algebra $A$, provided $M$ is a strongly $\mu$-free $A$-module.

## §3 A Step-Lemma and a Freeness-Proposition

**The $\mu$-$\aleph_0$-Step-Lemma 3.1:** *Let $H$ be a strongly $\mu$-free $R$-module of cardinality $\mu$ and $H^* = 0$. Moreover, let $F = \oplus_{n \in \omega} K_n$ be a direct sum of free modules $K_n$ of rank $\mu$. Choose a partial basis $x_{\alpha n}$ $(\alpha < \mu)$ for each $K_n$, let $K'_n = \langle x_{\alpha n} : \alpha < \mu \rangle$ and let $f \in K_0$ be any pure element with $K_0/fR$ free. Then we can find two extensions $F^i \supset F$ $(i = 0, 1)$ with the following properties.*

(1) *If $K_{\beta n} = \langle x_{\alpha i} : i < n$ or $\alpha < \beta \rangle = \oplus_{i < n} K'_i \oplus \langle x_{\alpha i} : i \geq n, \alpha < \beta \rangle$ then $F^i/K_{\beta n}$ is free $R$-module for each $n \in \omega, \beta < \mu$.*

(2) *$F^i/F \cong H$ $(i = 0, 1)$.*

(3) *If $\varphi \in F^*$ extends to both $\varphi^i \in F^{i*}$ $(i = 0, 1)$, then $f\varphi = 0$.*

**Remarks:** It follows from (1) that $F^i$ is free, because $K_{\beta n}$ is free as well. The particular case $\beta = 0$ with $K_{0n} = \oplus_{i < n} K_i$ looks more familiar in connection with older Step-Lemmas, but we need more – that is stated in (1). Conditions (1) and (2) say that $H$ has an $\langle \aleph_0 \rangle$–resolution for $\aleph_0$ as defined in (2.3).

Later on we will need a more general Step-Lemma as well, in order to deal with abelian groups of cardinality $\geq \aleph_{\omega+1}$. This however is postponed to Section 7.

**Proof.** Recall our representation (2.2) of $H$. We have

$$(i) \begin{cases} H_{\alpha+1} & = & \langle H_\alpha, t_{\alpha n} : n \in \omega \rangle \text{ with} \\ q_{\alpha n} t_{\alpha n+1} & = & t_{\alpha n} + \sum_{j < n_\alpha} t_{\alpha_j n_j} \text{ if } \alpha \in \boldsymbol{S} \text{ and} \\ H_{\alpha+1} & = & H_\alpha \oplus t_\alpha R \text{ if } \alpha \in \mu \backslash \boldsymbol{S} \end{cases}$$

where $q_{\alpha n} \in S$, $t_{\alpha_j n_j} \in H_{\alpha+1}$ for $\alpha_j < \alpha$. Let $C^i = \bigoplus_{\alpha \in \boldsymbol{S}} \bigoplus_{n \in \omega} s^i_{\alpha n} R \oplus \oplus_{\alpha \in \mu} \boldsymbol{S} s^i_\alpha R$ be freely generated by elements $s^i_{\alpha n}$ and $s^i_\alpha$ for $i = 0, 1$. Next we choose two submodules $N^i \subseteq F \oplus C^i$ by generators

$$N^0 = \langle q_{\alpha n} s^0_{\alpha n+1} - s^0_{\alpha n} - \sum_{j < n_\alpha} s^0_{\alpha_j n_j} - x_{\alpha n} : \alpha \in \boldsymbol{S}, \ n \in \omega \rangle$$

and

$$N^1 = \langle q_{\alpha n} s^1_{\alpha n+1} - s^1_{\alpha n} - \sum_{j < n_\alpha} s^1_{\alpha_j n_j} - x_{\alpha n} - f : \alpha \in \boldsymbol{S}, \ n \in \omega \rangle,$$

where $f \in K_0$ is taken from the hypothesis of the Lemma. If $F^i = F \oplus C^i/N^i$, then we identify $x \in F \oplus C^i$ with $x + N^i$ in $F^i$, hence $F \subset F^i$ $(i = 0, 1)$ and the following relations hold in $F^i$ by choice of $N^i$:

$$(r^i) \quad q_{\alpha n} s^i_{\alpha n+1} = s^i_{\alpha n} + \sum_{j < n_\alpha} s^i_{\alpha_j n_j} + x_{\alpha n} + f \cdot \delta_{1i} \qquad (i = 0, 1)$$

where $\delta_{11} = 1$ and $\delta_{10} = 0$.

In order to show (2), we consider the maps $\psi_i : F^i \to H$ defined by $\psi_i | F = 0$ and $s_{\alpha n}^i \psi_i = t_{\alpha n}$ for all $\alpha \in \boldsymbol{S}$, $n \in \omega$, $s_\alpha^i \psi_i = t_\alpha$ for $\alpha \in \mu \backslash \boldsymbol{S}$ and $i = 0, 1$. These maps preserve the relations $(r^i)$ and $(i)$, hence $\psi_i$ can be extended to an $R$-homomorphism $\psi_i : F^i \to H$. Its kernel is clearly $F$, $\psi$ is epic and (2) follows.

We have a decomposition into free summands $K_n' \oplus L_n = K_n$ where $K_n' = \langle x_{\alpha n} : \alpha < \mu \rangle$, and $L_n$ comes from the complement of the partial basis, hence

$(ii)$ $\begin{cases} F &=& \oplus_{m < n} K_m \oplus \oplus_{m \geq n} L_m \oplus \langle x_{\alpha m} : \alpha < \mu, \ m \geq n \rangle \text{ and} \\ F^0 &=& \oplus_{m < n} K_m \oplus \oplus_{m \geq n} L_m \oplus C \text{ where} \\ C &=& \langle x_{\alpha m}, s_{\alpha k}^0, \ s_\beta^0 : \alpha \in \boldsymbol{S}, \ \beta \in \mu \backslash \boldsymbol{S}, \ m \geq n, \ k \in \omega \rangle \end{cases}$

Using $(r^0)$ we see that the $x_{\alpha m}$'s can be expressed by $s_{\alpha k}^0$'s, hence $C$ is also generated by

$$C = \langle s_{\alpha k}^0, \ s_\beta^0 : \alpha \in \boldsymbol{S}, \ \beta \in \mu \backslash \boldsymbol{S}, \ k < \omega \rangle$$

If $i = 0$ is replaced by $i = 1$ we note that $f \in K_0$ and a similar submodule $C$ is a new complement satisfying (ii). In order to show (1), we consider the case $i = 0$ only. Let $\overline{\phantom{x}} : F^0 \to F^0 / K_{\beta n}$ denote the canonical epimorphism. The first step in proving that $\overline{F}^0$ is free is

$(iii)$ $\begin{cases} \overline{C}_{n\beta} := \langle \overline{s}_{\alpha k}^0, \ \overline{s}_{\alpha'}^0 : \ \alpha' \in \mu \backslash \boldsymbol{S} \text{ and } \alpha \in \boldsymbol{S} \text{ such that} \\ k < n \text{ or } k \geq n \text{ and } \alpha', \alpha < \beta \rangle \subseteq \overline{C} \text{ is freely generated.} \end{cases}$

Note that
$C_\gamma := \langle \overline{s}_{\alpha k}^0, \ \overline{s}_{\alpha'}^0 : \alpha \in \boldsymbol{S}, \ \alpha' \in \mu \backslash \boldsymbol{S} \text{ with } k < n, \ \alpha', \alpha < \gamma \text{ or } k \geq n, \ \alpha', \ \alpha < \beta \rangle$
$(\beta \leq \gamma < \mu)$ is an ascending, continuous chain of submodules of $\overline{C}_{n\beta}$ with $\overline{C}_{n\beta} = \bigcup_{\gamma < \mu} C_\gamma$. Because $|C_\gamma| \leq |\gamma| + |\beta| + \aleph_0 < \mu$ and $C_\gamma \subseteq \overline{F}^0$ is strongly $\mu$-free by (2), as an extension of a free module by a strongly $\mu$-free module, each $C_\gamma$ is free as well. Now (iii) will follow if we show

$$C_{\gamma+1}/C_\gamma \cong \langle \overline{s}_{\gamma n - 1}^0 + C_\gamma \rangle \quad (\cong R)$$

for each $\gamma \geq \beta$. Observe that this is trivial if $\gamma \notin \boldsymbol{S}$, hence we assume $\gamma \in \boldsymbol{S}$. For simplicity we will ignore the elements $s_\alpha^i$ $(\alpha \notin \boldsymbol{S})$! The relation $(r^0)$ under $\overline{\phantom{x}}$ become

$$q_{\gamma k} \overline{s}_{\gamma k + 1}^0 = \overline{s}_{\gamma k}^0 + \sum_{j < n_\gamma} \overline{s}_{\gamma j n_j}^0,$$

as in (i), and $C_{\gamma+1}/C_\gamma = \langle \overline{s}_{\gamma k}^0 + C_\gamma : \ k < n \rangle$. For large enough $q \in S$ we get $q \overline{s}_{\gamma n - 1}^0 = q' \overline{s}_{\gamma 0}^0 + \sum_j q_j \overline{s}_{\beta_j 0}$ and $\sum_j q_j \overline{s}_{\beta_j 0}$ is in $C_\gamma$, $q' | q$. Hence the $n$ generators above can be replaced by $C_{\gamma+1}/C_\gamma = \langle \overline{s}_{\gamma n - 1}^0 + C_\gamma \rangle$ and (iii) holds.
Next we show that

$(iv)$ $\begin{cases} \overline{C}_{n\beta} \oplus \overline{C}^{n\beta} = \overline{C} \text{ where } \overline{C}^{n\beta} = \langle \overline{s}_{\alpha k}^0 : \ \beta \leq \alpha < \mu, \ k \geq n \rangle \\ \text{is freely generated by the } \overline{s}_{\alpha k}^0 s. \end{cases}$

Clearly $\overline{C}$ is generated by $\overline{C}_{n\beta}$ and $\overline{C}^{n\beta}$.

Any linear combination $s = \sum r_{\alpha k} \overline{s}^0_{\alpha k} \in \overline{C}^{n\beta}$ may be multiplied by a large enough $q \in S$ such that the summand $q r_{\alpha k} \overline{s}^0_{\alpha k}$ by $(r^0)$ becomes a sum with summand $q' r_{\alpha k} x_{\alpha k}$. Now $s \in \overline{C}_{n\beta}$ is only possible if all the independent terms $q' r_{\alpha k} x_{\alpha k}$ ($\alpha \geq \beta$, $k \geq n$) vanish, hence $r_{\alpha k} = 0$ and $s = 0$. The sum (iv) must be direct. The freeness of $\overline{C}^{n\beta}$ follows similarly, again using independence of the $x_{\alpha k}$ under $^-$, and (iv) follows.

The quotient $\overline{F}^0 / \overline{C}_{n\beta} \cong \bigoplus_{m \geq n} L_m \oplus \overline{C}^{n\beta}$ is free by (iv) and $\overline{F}^0$ is free by (iii) and (1) is shown for $i = 0$.

In order to show (3) we consider any $\varphi \in F^*$ extending to $\varphi^i \in F^{i*}$ for $i = 0, 1$ and recall that $f \in F_0 \subseteq F$ is a prescribed basic element by hypothesis.

First we consider the pushout $D$ of $F \subset F^i$, which is

$$D = F^0 \oplus F^1 / \triangle \ \ with \ \triangle = \{(g, -g) : \ g \in F\} \cong F.$$

We also have a pushout extension $\varphi^i \subset \Phi \in D^*$ ($i = 0, 1$) by the pushout properties and $\varphi \subseteq \varphi^i$.

If $\overline{g} = (g, 0) + \triangle$, $\overline{g} = (0, g) + \triangle$ for $g \in F^0$ and $g \in F^1$ respectively, then we have canonical identifications $\overline{F^0} + \overline{F^1} = D$ with $\overline{F} = \overline{F^0} \cap \overline{F^1} = \{(g, 0) + \triangle : \ g \in F\} = \{(0, g) + \triangle : \ g \in F\}$. We omit "$^-$" in the following. Let $U = \langle F, d_{\alpha n}, d_{\alpha'} : \ \alpha \in \boldsymbol{S}, \ n < \omega, \ \alpha' \in \mu \backslash \boldsymbol{S} \rangle \subseteq D$ for $d_{\alpha n} = s^0_{\alpha n} - s^1_{\alpha n}$, $d_\alpha = s^0_\alpha - s^1_\alpha$. Observe that

$$U + F^i = F^0 + F^1 = D \ \text{and} \ U \cap F^i = F^0 \cap F^1 = F.$$

The desired relations from $(r^i)$ are

$$(v) \begin{cases} q_{\alpha n} d_{\alpha n+1} = d_{\alpha n} + \sum_{i < n_\alpha} d_{\alpha_i n_i} + f. \\ \text{The map } t_{\alpha n} \to d_{\alpha n} + \langle f \rangle \text{ induces } V / \langle f \rangle \cong H, \text{ with} \\ V = \langle f, \ d_{\alpha n}, \ d_{\alpha'} : \ \alpha \in \boldsymbol{S}, \ n \in \omega, \ \alpha' \in \mu \backslash \boldsymbol{S} \rangle. \end{cases}$$

In particular, $U/F \cong H$. Let $f\varphi = r \in R$ and consider the new homomorphism $\psi = \Phi - r : \ D \to R$.

Now we have $f\psi = f(\Phi - r) = 0$ and the relations (v) become

$$q_{\alpha n}(d_{\alpha n+1} \psi) = (d_{\alpha n} \psi) + \sum_{i < n_\alpha} (d_{\alpha_i n_i} \psi)$$

and $t_{\alpha n} \xrightarrow{\eta} d_{\alpha n} \psi$, $f\psi = 0$, induces an epimorphism

$$\eta : H \to V\psi$$

and $H^* = 0$ forces $V\psi = 0$. Hence $\Phi | V$ is multiplication by $r$, which by cardinality reasons is only possible for $r = 0$. We conclude $f\varphi = 0$.

The following Freeness-Proposition is a special case of a more general Freeness-Proposition. In order to make the idea of the construction more transparent, we postpone again the extension of (3.2) to Section 7. The more general Freeness-Proposition 7.2 involves a few more technical difficulties, due to more complicated notation. The Proposition will ensure that continuous chains of free submodules in our construction remain free at limit ordinals below $\mu^+$.

**The $\omega$-Freeness-Proposition 3.2:** *Let $\mu$ be a regular cardinal and $\delta < \mu^+$ be an ordinal with $\{G_i : i \leq \delta\}$ an ascending, continuous chain of submodules with the following properties for all $i < \delta$.*

(a) *$G_i$ is free of rank $\mu$.*

(b) *$G_{i+1}/G_i$ is strongly $\mu$-free and let $\boldsymbol{S} = \{i < \delta : G_{i+1}/G_i$ is not free $\}$*

(c) *There exist $G_i \subseteq H_i \subseteq G_{i+1}$ and $H_i^0$ a free summand of $G_i$ with free quotient such that $H_i = G_i \oplus_{H_i^0} H_i^1$ is a pushout with $H_i^0 \subseteq H_i^1$ free modules of rank $\mu$.*

(d) *$G_{i+1} = H_i \oplus K_i$ with $K_i$ freely generated by $x_{\alpha i}$ $(\alpha < \mu)$.*

(e) *If $i \notin \boldsymbol{S}$, then $G_i = H_i$ and $H_i^1 = H_i^0$ (hence $G_{i+1} = G_i \oplus K_i$).*

(f) *If $i \in \boldsymbol{S}$, then $\mathrm{cf}(i) = \omega$. Moreover*

    (i) *there exists $a_i \subseteq i$ of order type $\omega$ with $a_i \cap \boldsymbol{S} = \emptyset$ and $\sup a_i = i$.*

    (ii) *there exist $A_{ij} \subseteq \mu$ for all $j \in a_i$ such that $A_{ij} \cap A_{ij'} = \emptyset$ for $j \neq j'$ and $A_{ij} \cap A_{i'j'}$ is bounded in $\mu$ for $(ij) \neq (i'j')$.*

    (iii) *Let $K_j^\omega = \langle x_{\alpha j} : \alpha \in A_{ij} \rangle$, $H_i^0 = \bigoplus_{j \in a_i} K_j^\omega$. If $k \in a_i$, $\beta < \mu$ and $K_{\beta k} = \bigoplus_{j \in a_i}^{j < k} K_j^\omega \oplus \langle x_{\alpha j} : j \in a_i, \ j \geq k, \ \alpha < \beta \rangle$, then $H_i^1/K_{\beta k}$ is free.*

*It follows that $G_\delta$ is free as well.*

**Proof.** Enumerate $\boldsymbol{S} = \{i_\epsilon : \epsilon < \sigma\}$ for some $\sigma \leq \mu$. By induction on $\epsilon < \sigma$ choose ordinals $\alpha_\epsilon < \mu$ and suppose $\alpha_0, \cdots, \alpha_\nu$ are defined for all $\nu < \epsilon$. Define $\alpha_\epsilon$ subject to the following condition

    (+) For all $\nu < \epsilon$, $j \in a_{i_\nu} \cap a_{i_\epsilon}$ we have $A_{i_\nu j} \cap A_{i_\epsilon j} = \emptyset$.

This can be arranged: we collect all troublemakers at stage $\epsilon$, which is a set $M = \bigcup_{\nu < \epsilon} \bigcup_{j \in a_{i_\nu} \cap a_{i_\epsilon}} (A_{i_\nu j} \cap A_{i_\epsilon j}) \subseteq \mu$, and observe that $A_{i_\nu j} \cap A_{i_\epsilon j}$ is bounded below $\mu$ by (f)(ii), $\mu$ is regular and $|\epsilon| < \mu$. Hence $M$ is bounded by some $\alpha_\epsilon < \mu$ and (+) follows. Rename $\alpha_\epsilon = \alpha_{i_\epsilon}$. Then (+) can be restated as

    (+) For all $i \neq t \in \boldsymbol{S}$, $j \in a_i \cap a_t$ we have $A_{tj} \backslash \alpha_t \cap A_{ij} \backslash \alpha_i = \emptyset$.

If $i \in \mathbf{S}$, choose any fixed $\beta_i < \mu$ such that $\alpha_i < \beta_i$ and observe that $H_i^0 = K_i^* \oplus C_i$ by (f)(iii) for $k = min\ a_i$, $\beta = \beta_i$ and
$K_{k\beta} = K_i^* := \langle x_{\alpha j} :\ j \in a_i,\ \alpha < \beta_i \rangle$.
Now we are able to collect a basis of $G_\delta$: Choose all basis elements $x_{\alpha j}$ not involved in extension of $G_i$ at $i \in \mathbf{S}$, and we have the set

$$B = \{x_{\alpha j} :\ j < \delta,\ \text{there is no } i \in \mathbf{S} \text{ with } j \in a_i,\ \alpha \in A_{ij} \text{ and } \alpha > \beta_i\}.$$

Moreover $C = \bigoplus_{i \in \mathbf{S}} C_i$ is direct because of $(+)$, and each $C_i$ is free by (f)(iii). We may collect its basis $B_C$. Finally choose all basis elements $x_{\alpha j}$ $(\alpha < \beta_i)$ which are used for the distinct $K_i^*$, say $B'$. Then $B \cup B_C \cup B'$ is a basis for $G_\delta$.

# §4 The prediction principles

We will use prediction principles that hold under G.C.H. First recall some well-known definitions, tailored for our applications. Let $\kappa$ be a regular cardinal and $\theta = \{A_\nu :\ \nu < \kappa\}$ be a $\kappa$-filtration of $A = \bigcup_{\nu < \kappa} A_\nu$, i.e. $\theta$ is an ascending continuous chain with $A_0 = 0$ and $|A_\nu| < \kappa$ for all $\nu < \kappa$. Let $E \subseteq \kappa$. Then $\diamondsuit_\kappa(E)$ holds if there is a family $\{S_\nu \subseteq A_\nu :\ \nu \in E\}$ such that for all $X \subseteq A$ the set $\{\nu \in E :\ X \cap A_\nu = S_\nu\}$ is stationary in $\kappa$. Besides this diamond principle, we also consider the weak diamond $\Phi_\kappa(E)$. If $P_\nu : \mathcal{P}(A_\nu) \to 2 = \{0,1\}$ $(\nu \in E)$ is a given partition of subsets of $A_\nu$, then $\Phi_\kappa(E)$ provides a prediction function $\varphi : E \to 2$ such that for all $X \subseteq A$ the set $\{\nu \in E :\ P_\nu(X \cap A_\nu) = \varphi(\nu)\}$ is stationary in $\kappa$. In this case $E$ is called non-small. While $\diamondsuit_\kappa(E)$ holds in $V = L$ for all stationary subsets $E$ of regular, not weakly compact cardinals by Jensen [30], prediction principles are quite often valid under weaker assumptions, cf. Eklof, Mekler [15, pp.175-178], Shelah [35, p.376, Theorem 32]. See also Gregory [24].

**Proposition 4.1** (Shelah [35] ): *Suppose $\lambda = 2^\nu = \mu^+$ and for some regular cardinal $\kappa < \mu$ either*

(i) $\mu^\kappa = \mu$ *or*

(ii) *$\mu$ is singular with $cf\mu \neq \kappa$ and $|\delta|^\kappa < \mu$ for all $\delta < \mu$*

*Then $\diamondsuit_\lambda(E_\kappa)$ holds, where $E_\kappa = \{\alpha < \lambda :\ cf\alpha = \kappa\}$.*

In particular we have a

**Corollary 4.2:** *Assume G.C.H. and $\mu$ is some regular cardinal. Then $\diamondsuit_\mu E_\kappa$ holds for any regular cardinal $\kappa < \mu$.*

In the first application we will use $\Phi_\mu(E_{\aleph_0})$, a weak consequence of (4.2). Later we will replace $\kappa = \aleph_0$ by some regular $\kappa > \aleph_0$ (cf. §7). Moreover we will use

**Proposition 4.3:** *If $\lambda = \mu^+$ and $E \subseteq \lambda$ is non-small, then we can decompose $E = \coprod_{\beta < \kappa} E_\beta$ and each $E_\beta$ is non-small as well.*

**Proof.** See Eklof and Mekler [15, p.144].

# §5 Construction of modules with trivial dual

We apply the $\mu - \aleph_0$-Step-Lemma 3.1, the $\omega$-Freeness-Proposition 3.2 and combinatorial ideas from §4. The arguments are however partly standard (for working in $L$), so we keep them short. The reader may consult [15], [38] or [8] for details.

Let $\lambda = \mu^+ = 2^\mu$ for some regular cardinal $\mu$, and let $E = \{\alpha \in \lambda,\ cf\,\alpha = \omega\}$. Recall that $E$ is non-small by (4.1), for instance. Hence we assume $\Phi_\lambda(E)$, the weak diamond and decompose $E = \coprod_{\beta < \lambda} E_\beta$ into non-small sets $E_\beta$, cf. (4.3). For each $i \in E$ we also fix $a_i = \{i_n :\ n \in \omega\}$, a strictly increasing sequence $i_n < i$ of successor ordinals $i_n$ with $sup\,a_i = i$, hence $E \cap a_i = \emptyset$. For each $j \in a_i$ we find $A_{ij} \subseteq \mu$ with:

(a) $A_{ij}$ *unbounded in $\mu$ and $A_{ij} \cap A_{i'j'}$ bounded in $\mu$ if $(ij) \neq (i'j')$ and*

(b) $A_{ij} \cap A_{ij'} = \emptyset$ *for $j \neq j' \in a_i$.*

Choose a set of $2^\mu = \lambda$ subsets $X', X \subseteq \mu$ with $|X \cap X'| < \mu$ for distinct pairs $X, X'$. Decompose each $X$ into $\omega$ subsets of cardinal $\mu$ to define the $A_{ij}$'s. Hence (b) follows and (a) holds because $\mu$ is regular. Suppose now that $H$ is a given strongly $\mu$-free $R$-module with $|H| = \mu$ and trivial dual $H^* = Hom\,(H, R) = 0$.

We want to derive the existence of a strongly $\lambda$-free $R$-module $G$ of cardinality $\lambda$ such that $G^* = 0$.

We have to fix the "prediction" and choose a $\lambda$-filtration $\{G_i :\ i < \lambda\}$ such that $G_0 = 0$ and $|G_i| = |G_{i+1} \backslash G_i| = \mu$ for all $i < \lambda$. Let $\{g_\nu :\ \nu < \lambda\} = G$ be an enumeration of $G = \bigcup_{i < \lambda} G_i$ such that $g_\beta \in G_i$ for all $i \in E_\beta$. We define partitions $P_i^\beta : \mathcal{P}(G_i \times R) \to 2$ for all $\beta < i$, $i \in E_\beta$ predicting homomorphisms in $G_i^*$.

Formally we also must predict an $R$-module structure on $G_i$, hence $G_i \times R$ must be replaced by $G_i \times G_i \times \cdots \times R$ to take care of this; however we will ignore this for simplicity; compare [8] or [19, p.284].

If $X \subseteq G_i \times R$, then define $P_i^\beta(X) = 0$ if the following holds.

(1) $\left\{\begin{array}{l} G_i, G_j, G_i/G_j \text{ are free } R\text{-modules for } j \in a_i \text{ and } g_\beta \in G_{j_0}, \text{ with } G_{j_0}/g_\beta R \\ \text{free, say } j_0 = \min a_i. \text{ Moreover } G_{j+1} = G_j \oplus K_j \text{ with } K_j \\ \text{freely generated by some } x_{\alpha j} \ (\alpha < \mu). \end{array}\right.$

(2) $\left\{\begin{array}{l} X \text{ is (as a graph) a homomorphism } h_X \in G_i^*. \text{ If } K_j^\omega = \langle x_{\alpha j} : \alpha \in A_{ij} \rangle \\ (j \in a_i) \text{ and if we identify } (G_i, \oplus_{j<n} K_j^\omega, g_\beta, h_X : n \in \omega) \text{ with} \\ (F, \oplus_{j<n} K_j', f, \varphi : n \in \omega) \text{ in (3.1), then we} \\ \textbf{require} \text{ that } \varphi \text{ does not extend to } F^{0*}. \end{array}\right.$

From $\Phi_\lambda(E_\beta)$ we have prediction functions $\varphi_\beta : E_\beta \to 2$ such that

$$\mathcal{X}_\beta(X) := \{i \in E_\beta, \ P_i^\beta(X \cap G_i \times R) = \varphi_\beta(i)\}$$

is stationary in $\lambda$ for all $X \subseteq G$, $\beta < \lambda$. The structure of the $R$-module $G = \bigcup_{i<\lambda} G_i$ is given inductively on the continuous chain $\{G_i : i < \lambda\}$ subject to the following conditions.

(i) $G_i$ is a free $R$-module for any $i \in \lambda$.

(ii) $\left\{\begin{array}{l} G_{i+1} = H_i \oplus K_i \text{ where } K_i \text{ is freely generated by elements} \\ x_{\alpha i} \ (\alpha < \mu), \ G_i \subseteq H_i \subseteq G_{i+1}, \ G_{i+1}/G_i \text{ is strongly } \mu - \text{free} \\ \text{and } H_i = G_i \oplus_{H_i^0} H_i^1 \text{ is a pushout for some free summand} \\ H_i^0 = G_i \cap H_i^1 \text{ of } G_i \text{ with free quotient and } H_i^0 \subseteq H_i^1 \text{ all} \\ \text{free of rank } \mu. \end{array}\right.$

(iii) $\left\{\begin{array}{l} \text{If } i \in E_\beta \text{ and (1) holds, then we choose } H_i \text{ as an extension} \\ F^{\varphi_\beta(i)} \text{ given by the identification } (F, \oplus_{j<n} K_j', f, \varphi : n \in \omega) \\ \text{in (3.1) with } (G_i, \oplus_{j<n} K_j^\omega, g_\beta, h_X : n \in \omega). \end{array}\right.$

Let $H_i = G_i$ otherwise. Hence $G$ is defined.

Note that $G_\delta = \bigcup_{i<\delta} G_i$ is a free $R$-module for limit ordinals $\delta < \lambda$ as follows from (3.2). We conclude that $G$ is a strongly $\lambda$-free $R$-module of cardinality $\lambda$ by (3.1) and the construction above.

If $\psi \in G^*$ and $\psi \neq 0$, then $C = \{i < \lambda : \psi|G_i \neq 0\}$ is a cub. Moreover $f\psi \neq 0$ for some $f \in G$ with $G_i/fR$ free for large enough $i \in C$, and $f = g_\beta \in G_i$ for some $\beta < \lambda$. Choose an ordinal $i \in C \cap \mathcal{X}_\beta(\psi)$ and observe that $\psi|G_i$ extends to $G_{i+1}$, hence $P_i^\beta(\psi|G_i) = 1$ and $G_{i+1} = F^1$ as in the Step-Lemma 3.1. However $P_i^\beta(\psi|G_i) = 1$ and (2) tell us that $\psi|G_i$ extends to $F^0$ as in (3.1) as well. Therefore $f\psi = 0$ by (3.1), a contradiction. Hence $G^* = 0$ follows.

We conclude the following results.

**Theorem 5.1:** *Let $R$ be a countable commutative ring as in §2. If $H$ is a strongly $\mu$-free $R$-module of regular cardinality $\mu$ with trivial dual $H^* = 0$ and $\lambda = 2^\mu = \mu^+$, then there exists a strongly $\lambda$-free $R$-module of cardinality $\lambda$ with trivial dual.*

**Corollary 5.2:** *If $2^{\aleph_i} = \aleph_{i+1}$ for $i = 0, \cdots, n-1$, then there exists a strongly $\aleph_n$-free abelian group of cardinality $\aleph_n$ with trivial dual.*

The corollary follows by induction. If $i = 0$, then any subgroup of $\mathbb{Q}$ different from $\mathbb{Z}$ serves as $H$ in (5.1). From (5.1) we obtain a strongly $\aleph_1$-free abelian group of cardinality $\aleph_1$ with trivial dual, and so on.

## §6 Endo-rigid modules

The following result establishes the existence of $R$-modules with endomorphism ring $R$, the so-called rigid $R$-modules. Such modules are known in ZFC, see e.g. Corner [2] or [20], [21], [40]. However we will also assume that the $R$-modules in question are almost free. A result in $ZFC + \mathrm{MA}$ (= Martin's axiom) shows that extra axioms are needed, because such modules turn out to be separable under ZFC + MA, see Göbel, Shelah [22]. Assuming $V = L$, almost free modules were constructed in Dugas [7], Shelah [38] and Dugas, Göbel [8] several years ago. Here we have to work harder to derive the existence under ZFC and a weak form of G.C.H.

**Theorem 6.1:** *Let $R$ be a countable, commutative ring as in (2.1) and assume $\lambda = 2^\mu = \mu^+$ for some regular cardinal $\mu$. Suppose there exists a strongly $\mu$-free $R$-module $H$ of cardinality $\mu$ with $\mathrm{End}\,H = R$. Then we can find a strongly $\lambda$-free $R$-module $G$ of cardinality $\lambda$ with $\mathrm{End}\,G = R$.*

The proof is – as in (5.1) – divided into three parts, a Step-Lemma similar to (3.1), however with crucial, not-so-obvious changes. The construction of $G$ and the actual proof of the claims in (6.1) are similar to §5. We will keep the proof short. However, the reader familiar with §§3-5 will follow the arguments easily. We begin with the Step-Lemma.

**The $\mu$-$\aleph_0$-Step-Lemma 6.2:** *Let $R$ be as in §2 and let $H$ be a strongly $\mu$-free $R$-module of cardinality $\mu$ with $\mathrm{End}\,H = R$. Moreover, let $F = \oplus_{n \in \omega} K_n$ be a direct sum of free modules $K_n$ of rank $\mu$. Choose a partial basis $x_{\alpha n}$ $(\alpha < \mu)$ for each $K_n$, let $f \in K_0$ and $h \in H$ such that $K_0/fR$ is free and $H/hR$ is strongly* *μ-free. Then we can find two extensions $F \subseteq F^i$ $(i = 0, 1)$ with the following properties.*

(1) *If $K_{\beta n} = \langle x_{\alpha i} : i < n \text{ or } \alpha < \beta \rangle$, then $F^i / K_{\beta n}$ is a free $R$-module for each $n \in \omega$, $\beta < \mu$.*

(2) *$F^i / F \cong H/hR$ is strongly $\mu$-free.*

(3) *If $\varphi \in End\, F$ extends to both $\varphi^i \in End\, F^i$ $(i = 0, 1)$, then $f\varphi \in fR$.*

**Proof.** Let $H' = H/hR$ and choose a $\mu$-filtration $H = \bigcup_{i < \mu} H_i$ for strong $\mu$-freeness with $h \in H_0$. Since $H/hR$ is strongly $\mu$-free, also $H'_i = H_i/hR$ is free and $\{H'_i : i < \mu\}$ constitutes a $\mu$-filtration of $H'$ for strong $\mu$-freeness.

Let $\boldsymbol{S} = \{\alpha < \mu \text{ with } H_{\alpha+1}/H_\alpha \text{ not free }\}$. Next we apply the representation (2.2) to $H'$ and obtain

$$(*)\ \begin{cases} \text{a continuous chain of free submodules } H_\alpha\ (\alpha < \mu) \text{ of} \\ H = \bigcup_{\alpha < \mu} H_\alpha \text{ such that } H_0 = hR \text{ and } H_{\alpha+1} = H_\alpha \oplus t_\alpha R \\ \text{if } \alpha \notin \boldsymbol{S} \text{ or } H_{\alpha+1} = \langle H_\alpha, t_{\alpha n} : n < \omega \rangle \text{ with} \\ q_{\alpha n} t_{\alpha n+1} = t_{\alpha n} + \sum_{j < n_\alpha} t_{\alpha_j n_j} + r_{\alpha n} h \text{ for} \\ q_{\alpha n}, r_{\alpha n} \in S,\ t_{\alpha_j n_j} \in H_{\alpha_j + 1} \text{ with } \alpha_j < \alpha,\ n < \omega. \end{cases}$$

Similar to (3.1), but identifying $f$ and $h$, we define two extensions $F^i$ of $F$ for $i = 0, 1$. Let $C^i = \bigoplus_{\alpha \in \boldsymbol{S}} \bigoplus_{n < \omega} s_{\alpha n}^i R \oplus \bigoplus_{\alpha \in \mu \backslash \boldsymbol{S}} s_\alpha^i R$ be freely generated and $N^i \subseteq F \oplus C^i$ such that

$$N^0 = \langle q_{\alpha n} s_{\alpha n+1}^0 - s_{\alpha n}^0 - \textstyle\sum_{j < n_\alpha} s_{\alpha_j n_j}^0 - x_{\alpha n} : \alpha \in \boldsymbol{S},\ n \in \omega \rangle$$
$$N^1 = \langle q_{\alpha n} s_{\alpha n+1}^1 - s_{\alpha n}^1 - \textstyle\sum_{j < n_\alpha} s_{\alpha_j n_j}^1 - x_{\alpha n} - r_{\alpha n} f : \alpha \in \boldsymbol{S},\ n \in \omega \rangle$$

Let $F^i = F \oplus C^i / N^i$ and identify $x \in F \oplus C^i$ with $x + N^i$ in $F^i$. The new relations are

$$(r^i)\quad q_{\alpha n} s_{\alpha n+1}^i = s_{\alpha n}^i + \sum_{j < n_\alpha} s_{\alpha_j n_j}^i + x_{\alpha n} - r_{\alpha n} \delta_{1i} f \quad (i = 0, 1)$$

The maps $\psi_i : F^i \to H/hR$ given by $\psi_i | F = 0$, $s_\alpha^i \psi_i = t_\alpha + hR$ and $s_{\alpha n}^i \psi_i = t_{\alpha n} + hR$ preserve the relation in $(*)$ modulo $hR$ and are epic, hence (2) follows.

The proof of (1) is the same as in (3.1).

Suppose $\varphi \in End\, F$ extends to $\varphi^i \in End\, F^i$ $(i = 0, 1)$. Again we consider the pushout $F^0 + F^1 = D$ with $F = F^0 \cap F^1$ as in (3.1), and let
$U = \langle F, d_{\alpha n}, d_{\alpha'} : \alpha \in \boldsymbol{S},\ n \in \omega,\ \alpha' \in \mu \backslash \boldsymbol{S} \rangle \subseteq D$
for $d_{\alpha n} = s_{\alpha n}^0 - s_{\alpha n}^1$ and $d_\alpha = s_\alpha^0 - s_\alpha'$. We also have $U + F^i = D$ and $U \cap F^i = F$, however the relations for $d_{\alpha n}$ from $(r^i)$ are
$q_{\alpha n} d_{\alpha n+1} = d_{\alpha n} + \sum_{j < n_\alpha} d_{\alpha_j n_j} + r_{\alpha n} f$. Hence

$$V = \langle f, d_{\alpha n}, d_{\alpha'} : \alpha \in \boldsymbol{S},\ n \in \omega \text{ and } \alpha' \in \mu \backslash \boldsymbol{S} \rangle \cong H$$

by $(*)$. We find a projection $\pi : D \to V$, because - as in (3.1) - $fR$ has a complement $L_{-1}$ in $K_0 \subseteq F$ (cf. 3.1) and we find complements $L_n$ of each $K_n$ in $K_{n+1}$, hence $D \subseteq V \bigoplus \bigoplus_{n \geq -1} L_n$ gives $\pi$. The pushout $D$ provides $\Phi \in End\ D$ extending $\varphi^i$ $(i = 0, 1)$ and $\Phi' = \Phi | V \cdot \pi$ is an endomorphism of $V$. Finally $End\ H = R$ forces $f\varphi = f\Phi \in fR$.

**Proof (6.1):** The proof is similar to §5, however (3.1) is replaced by (6.2). We only indicate the changes. Assuming $\lambda = \mu^+ = 2^\mu$ we prepare the application of $\Phi_\lambda(E)$ for $E = \{\alpha \in \lambda, \ cf\alpha = \omega\}$. Condition (2) in §5 must be replaced by a new condition

(2) $\begin{cases} X \text{ is an endomorphism } h_X \in End\ G_i.\ \text{If } K_j^\omega = \langle x_{\alpha j} : \ \alpha \in A_{ij} \rangle \\ (j \in a_i)\ \text{and if we identify } (G_i,\ \oplus_{j<n}K_j^\omega,\ g_\beta,\ h_X :\ n \in \omega)\ \text{with} \\ (F,\ \oplus_{j<n}K_j',\ f,\ \varphi :\ n \in \omega)\ \text{in}\ (6.2),\ \text{then we} \\ \textbf{\textit{require}}\ \text{that } \varphi \text{ does not extend to } End\ F^0. \end{cases}$

In the construction of $G$ we must change condition (iii). Now we have

(iii) $\begin{cases} \text{If } i \in E_\beta \text{and (1) holds, then we choose } H_i \text{as an extension} \\ F^{\varphi_\beta(i)} \text{given by the identification } (F,\ \oplus_{j<n}K_j',\ f,\ \varphi :\ n \in \omega) \\ \text{in (6.2) with } (G_i, \oplus_{j<n}K_j^\omega,\ g_\beta,\ h_X :\ n \in \omega). \end{cases}$

As in §5, the final $R$-module $G$ has size $\lambda$ and is strongly $\lambda$-free. In particular $R \subseteq End\ G$ by scalar multiplication. If $\psi \in End\ G \backslash R$, then we can find a cub

$$C = \{i < \lambda :\ \psi | G_i = \psi_i \in End\ G_i \backslash R\}.$$

Moreover we can choose $f \in G_i$ with $G_i/fR$ free such that $f\psi \notin fR$ and $f = g_\beta$ by enumeration. Choose an ordinal $i \in C \cap \mathcal{X}_\beta(\psi)$ and observe that $\psi_i$ extends to an endomorphism $\psi' \in End\ G_{i+1}$. The latter follows from the fact that $G_{i+1}\psi \subseteq G_j$ for some $j < \lambda$ and $G_{i+1}$ is a summand of $G_j$, so $\psi' = \psi | G_{i+1} \cdot \pi$ for some projection $\pi$. Hence $P_i^\beta(\psi_i) = 1$ and $G_{i+1} = F^1$ as in the Step-Lemma 6.2. However $P_i^\beta(\psi_i) = 1$ and (2) tell us that $\psi_i$ extends to $F^0$ as in (6.2) as well. Therefore $f\psi \in fR$ by (6.1), a contradiction, and $End\ G = R$ follows.

Inductively we derive a

**Corollary 6.3:** *If* $2^{\aleph_i} = \aleph_{i+1}$ *for* $i = 0, \cdots, n-1$,
*then there exists a strongly* $\aleph_n$-*free abelian group* $G$ *with* $|G| = \aleph_n$ *and* $End\ G = R$.

We also have a corollary from the proof of (6.1) observing that the $R$-algebra $R$ can be replaced by an $R$-algebra $A$.

**Theorem 6.1\*:** *Let $R$ be a countable, commutative ring which is $S$-reduced and $S$-torsion-free for some multiplicatively closed set $S$ of regular elements in $R$. Assume $\lambda = 2^\mu = \mu^+$ for some regular cardinal $\mu$ and suppose there exists a strongly $\mu$-free $A$-module $H$ of cardinality $\mu$ with $End_R H = A$ and $|A| < \mu$, $A_R$ free. Then we can find a strongly $\lambda$-free $A$-module $G$ of cardinality $\lambda$ with $End_R G = A$.*

Hence we have a corollary which ensures counterexamples for Kaplansky's test problems, cf. Corner [3] and Corner, Göbel [4]. The next corollary strengthens (6.3).

**Corollary 6.3\*:** *Let $R$ be as in (6.1\*) and assume $2^{\aleph_i} = \aleph_{i+1}$ for $i = 0, \cdots, n-1$. If $A$ is a countable, free $R$-algebra, then there exists a strongly $\aleph_n$-free $R$-module $G$ with $|G| = \aleph_n$ and $End\, G = A$.*

# §7 Almost free abelian groups which are of cardinality $> \aleph_\omega$

In order to derive results like the existence of strongly $|G|$-free groups $G$ with $|G| = \aleph_{\omega+1}$ and $End\, G = \mathbb{Z}$ under G.C.H., the induction used in §§3, 4 and 5 breaks down at $\aleph_\omega$ and we must extend the basic results of the earlier sections in order to overcome this difficulty. The changes will be sketched below. If $\mu$ is a regular cardinal, then let $\boldsymbol{r}(\mu)$ be the set of all decreasing sequences $u = \langle \mu_1, \cdots, \mu_n \rangle$ of regular cardinals $\mu \geq \mu_1 > \cdots > \mu_n$. Recall the notion of a $u$-resolution from (2.3) and of the freeness type of an $R$-module $H$ from §2. In case $End\, H = A$ we assume that $H$ is a strongly $\mu$-free $A$-module and replace $R$ by $A$ in the above definitions. If $H_0/hA$ is $A$-free for $H_0$ as in (2.2) and $h \in H_0$, then $ftH = ft(H/hA)$ follows. Observe that we showed in (3.1) (1) and (2) that a strongly $\mu$-free $R$-module $H$ of freeness type $u = \langle \mu, \aleph_0 \rangle$ has a $(\aleph_0)$-resolution. This implication does not depend on additional properties of $H$ as $H^* = 0$ or $End\, H = A$.

Moreover, we must refine the use of $\omega$-sequences converging to particular ordinals. If $u = \langle \mu_1, \cdots, \mu_n \rangle \in \boldsymbol{r}(\mu)$, then let $\times u = \mu_1 \times \cdots \times \mu_n$ be equipped with the lexicographical order. We consider strictly order preserving monomorphisms $\alpha : \times u \to \mu$ and let $a = Im(\alpha)$ and $i = sup\, a$ which is cofinal to $\mu_1$.
We say that $X \subseteq a$ is bounded (with respect to $\alpha$) if there is $\beta = (\beta_1, \cdots, \beta_n) \in \times u$ such that $X \subseteq Im(\alpha|\beta)$.
We may view $\alpha$ as an "$n$-ladder" at $a$. A bounded set $X$ is bounded "on each ladder" in an obvious sense. If $n = 1$, then $a$ is a $\mu_1$-sequence representing $cf(i) = \mu_1$. The bounded subsets of $a$ constitute an ideal of the set of subsets $\mathcal{P}(a)$ of $a$. The set $a \subseteq i$ will be used to find an extension on $G_i \subseteq G_{i+1}$ which will lead to the desired module $G = \bigcup_{i<\lambda} G_i$. If $u = \langle \aleph_0 \rangle$, as is the case considered in §3-6, then the 1-ladder system is an $\omega$-sequence $a_i$ at each ordinal $i$ with $cf(i) = \omega$ which has an additional property (+) (§3) which ensures freeness of $G_\delta$ at limit ordinals $\delta < \mu^+$. A similar condition is required for a general $u$:

$(+)$ $\begin{cases} \text{If } cf(i) = \mu_1, \text{ then there exists an } n\text{-ladder system } a_i \subseteq i \\ \quad \text{such that for each } \delta < \mu^+ \text{ there exist bounded subsets } b_i \subseteq a_i \\ \quad \text{(as defined above) such that the sets } (a_i \backslash b_i) \ i \in S_{\mu_1}, \ i < \delta \\ \quad \text{are pairwise disjoint.} \end{cases}$

This is ensured by the following combinatorial result of Shelah, which follows from the main theorem of [35].

**Theorem 7.1** (ZFC + G.C.H.): *If $S = \{i \in \aleph_{\omega+1}, \ cf(i) = \omega_1\}$ then $\diamondsuit_{\aleph_{\omega+1}}(S)$ holds. Moreover, for any $i \in S$ we can find a 2-ladder system $\alpha^i : \omega_1 \times \omega \to i$ such that the sets $a_i = Im(\alpha^i)$ satisfy (+).*

The $\mu$-$\aleph_0$-Step-Lemma 3.1 and (6.2) can be refined to

**The $\mu$-$u$-Step-Lemma 7.2:** *Let $(R, S)$ be as in (2.1) and $A$ be a countable, free $R$-algebra. We will distinguish two cases (a) and (b). In case (a) let $H = H'$ be a strongly $\mu$-free $R$-module of cardinality $\mu$ of freeness type $ftH = u$. In case (b) let $H$ be a strongly $\mu$-free $A$-module of cardinality $\mu$ with $End_R H = A$ and $h \in H_0$ (in (2.2)) such that $ftH = ftH' = u$ for $H' = H/hA$.*
*Then we can find two $u$-resolutions of $H'$ for $\mu_1$ (with respect to $R$ in case (a) and to $A$ in case (b)) of the form $F^j/K$ $(j = 0, 1)$ such that for any basic element $f \in K_0$ the following holds*

(a) *If $H^* = 0$ and $\varphi \in K^*$ extends to both $\varphi^j \in (F^j)^*$ $(j = 0, 1)$, then $f\varphi = 0$.*

(b) *If $\varphi \in End_R(K)$ extends to both $\varphi^j \in EndF^j$ $(j = 0, 1)$, then $f\varphi \in fA$.*

**Proof.** We consider case (a) only and must find suitable extensions $F^j$ of $K$. Recall that $K = \bigoplus_{i<\mu} K_i$ for $\mu = \mu_1$ is a direct sum of free $R$-modules $K_i$ of rank $\mu$, cf. (2.3). We will define $F^j$ for $j = 0, 1$. We fix an $n$-ladder $a \subseteq \mu$ and write $K$ as a continuous increasing chain $K = \bigcup_{i<\mu} F_i$ of free $R$-modules of rank $\mu$ such that $F_{i+1} = F_i \oplus K_i$ for all $i \in a$. The module $H$ can also be represented by (2.2) using $a$ in the form $H = \bigcup_{\alpha \in a} H_\alpha$. Using the ordering on $a$ (induced by $\mu$ and also by $\times u$), equations (2.2) can be expressed as

$$q_\alpha t_{\alpha+1} = t_\alpha + \Sigma_{\alpha_j < \alpha} t_{\alpha_j}.$$

The free resolution is now similar to $(3.1)(r^j)$. Defining $C^j$ and $N^j$ according to (3.1) we let $F^j$ (j = 0, 1) be generated by $K$ and new elements $s^j_\alpha$ $(\alpha \in a)$ subject to the relations

$$q_\alpha s^j_{\alpha+1} = s^j_\alpha + \Sigma_{\alpha_j < \alpha} s^j_{\alpha_j} + x_\alpha + \delta_{j1} f,$$

where $x_\alpha \in K_\alpha$ are fixed basic elements. Arguments as in (3.1) will give (a) and (b).

In order to ensure that unions at limit ordinals remain free in the construction of §5, the $\omega$-Freeness-Proposition 3.2 was used. At limit ordinals not cofinal to $\omega$, the construction was trivial. This is different here. If $u = \langle \mu_1, \cdots, \mu_n \rangle$ is a certain type of freeness, then the new construction will be non-trivial at limit ordinals cofinal to some $\mu_i$. Hence we have to work at those limit ordinals to ensure freeness of the constructed $\mu^+$ filtration. We need the generalized

**$u$-Freeness-Proposition 7.3:** *Let $\mu$ be a regular cardinal,*
*$u = \langle \mu_1, \cdots, \mu_n \rangle \in \boldsymbol{r}(\mu)$ some freeness type with $\mu_1 < \mu$ and $\delta < \mu^+$ an ordinal with $\{G_i :\ i \leq \delta\}$ an ascending, continuous chain of submodules with the following properties for all $i < \delta$.*

(a) $G_i$ *is free of rank $\mu$.*

(b) $G_{i+1}/G_i$ *is strongly $\mu$-free and let* $\boldsymbol{S} = \{i < \delta :\ G_{i+1}/G_i$ *is not free* $\}$.

(c) *There exist $G_i \subseteq H_i \subseteq G_{i+1}$ and $H_i^0$ a free summand of $G_i$ with free quotient such that $H_i = G_i \oplus_{H_i^0} H_i^1$ is a pushout with $H_i^0 \subseteq H_i^1$ free modules of rank $\mu$.*

(d) $G_{i+1} = H_i \oplus \bigoplus_{\tau \in u} K_i^\tau$ *with $K_i^\tau$ freely generated by $x_{\alpha i}^\tau$ $(\alpha < \mu)$.*

(e) *If $i \in \boldsymbol{S}$, then $G_i = H_i$ and $H_i^1 = H_i^0$ (hence $G_{i+1} = G_i \oplus \bigoplus_{\tau \in \boldsymbol{u}} K_i^\tau$).*

(f) *If $i \notin \boldsymbol{S}$, then $cf(i) = \kappa \in u$. Moreover*

    (i) *there exists $a_i \subseteq i$ of order type $\kappa$ with $a_i \cap \boldsymbol{S} = \emptyset$ and sup $a_i = i$.*

    (ii) *there exist $A_{ij} \subseteq \mu$ for all $j \in a_i$ such that $A_{ij} \cap A_{ij'} = \emptyset$ for $j \neq j'$ and $A_{ij} \cap A_{i'j'}$ is bounded in $\mu$ for $(ij) \neq (i'j')$.*

    (iii) *Let $K_j^\kappa = \langle x_{\alpha j}^\kappa :\ \alpha \in A_{ij} \rangle$, $H_i^0 = \bigoplus_{j \in a_i} K_j^\kappa$. If $k \in a_i$, $\beta < \mu$ and $K_{\beta k}^\kappa = \bigoplus_{j \in a_i}^{j < k} K_j^\kappa \oplus \langle x_{\alpha j}^\kappa :\ j \in a_i,\ j \geq k,\ \alpha < \beta \rangle$, then $H_i^1/K_{\beta k}^\kappa$ is free.*

*It follows $G_\delta$ is free as well.*

The proof of (7.3) follows the proof of (3.2) up to minor changes.

Combining (7.1) and our next result (7.5) we can pass $\aleph_\omega$ and obtain the existence of certain modules of cardinal $\lambda = \aleph_{\omega+1}$.

**Corollary 7.4** (ZFC + G.C.H.): *If $(R, S)$ is as in §2 and $A$ is a countable, free $R$-algebra, then we can find strongly $\aleph_{\omega+1}$-free $R$-modules $G$ of cardinality $\aleph_{\omega+1}$ with $End_R G = A$ and of freeness type $\langle \aleph_{\omega+1}, \aleph_1 \rangle$.*

We would like to remark that a further extension of (7.5) would provide modules as in (7.4) for cardinals $\aleph_{\omega k+1}$ $(k \in \omega)$. Clearly (7.4) in conjunction with (5.1) and (6.1*), respectively, provides examples for all $\aleph_{\omega+n}$ $(n \geq 1)$. A similar result then follows for all $\aleph_{\omega k+n}$ $(k, n \geq 1)$. Magidor and Shelah [31] show, however, that there are models of ZFC + G.C.H. in which $\aleph_{\omega^2+1}$-free abelian groups of cardinality $\aleph_{\omega^2+1}$ are free. This shows that the results above cannot be extended any further.

**Theorem 7.5** (ZFC + G.C.H.): *Suppose that the following conditions hold.*

(1) *Let $(R, S)$ be as in §2 and $A$ a countable, free $R$-algebra.*

(2) *There exists a strongly $\mu$-free $A$-module $H$ of cardinality $\mu$ with $End_R H = A$ of freeness type $u = \langle \mu_1, \cdots, \mu_n \rangle$.*

(3) *Let $\lambda$ be a regular, uncountable cardinal with a stationary, non-small subset $\boldsymbol{S}$ of ordinals cofinal to some cardinal $\mu_1$.*

(4) *There exists an $n$-ladder system $a_i \subseteq i$ $(i \in \boldsymbol{S})$ with $(+)$ and $\sup a_i = i$, $a_i \cap \boldsymbol{S} = \emptyset$.*

*If $\lambda = \mu^+$, then there exists a strongly $\lambda$-free $R$-module $G$ with $|G| = \lambda$, $ft(G) = \lambda^\wedge u$ and $End_R G = A$.*

**Remark:** A similar result holds for modules with trivial dual.

**Proof.** We indicate the changes in the construction of $G$, which are based on (2.2) using $u^n = \times u$ and the $n$–ladder $\alpha^i : \times u \to i$ with $Im(\alpha^i) = a_i$. Recall that $cf(i) = \mu$. From $a_i \cap \boldsymbol{S} = \emptyset$ and the inductive construction of $F = \cup_{j < \lambda} F_j$ we have $F_{\alpha+1} = F_\alpha \oplus K_\alpha$ for $\alpha \in a_i$. We choose basic elements $x_\alpha \in K_\alpha$ (with respect to $A$ in case of $End_R G = A$ and with respect to $R$ for $G^* = 0$) which are used in (7.2) in order to find $F_{i+1}^j$ $(j = 0, 1)$ with the desired properties. Moreover, (7.3) takes care of freeness at limit ordinals. The prediction principle (7.5)(c) picks the right branch in $2^\lambda$ which leads to $End_R G = A$ and $G^* = 0$, respectively.